\documentclass[12pt]{article}
\usepackage[final]{epsfig}
\usepackage{graphics}
\usepackage{amsmath}
\usepackage{amsfonts}
\usepackage{latexsym}
\usepackage{amssymb}
\usepackage{graphicx}
\usepackage{url}
\usepackage{epstopdf}
\usepackage{hyperref}
\usepackage{color}
\usepackage{marginnote}
\usepackage{a4wide}

\newtheorem{lemma}{Lemma}

\newtheorem{remark}[lemma]{Remark}

\newtheorem{theorem}{Theorem}

\newcommand{\g}{{\gamma}}

\newcommand{\proofend}{$\Box$\bigskip}

\def\proof{\paragraph{Proof.}}

\begin{document}

\title{A 4-point theorem: still another variation on an old theme}

\author{
Serge Tabachnikov
 \footnote{
Department of Mathematics,
Pennsylvania State University,
University Park, PA 16802,
USA;
tabachni@math.psu.edu}
}

\date{\today}

\maketitle

The subject of this article belongs to a ``neighborhood" of the 4-vertex theorem that, in the simplest form, states that the curvature of a plane oval (a smooth closed  curve with  positive curvature) has at least four critical points. Since its publication by S. Mukhopadhyaya in 1909, this result and its ramifications have generated a vast literature. We give but one reference, Lecture 10 of \cite{FT}. 

In what follows, we freely use basic facts of elementary differential geometry of the sphere and the hyperbolic plane; we omit references to numerous textbooks on the subject.

Our starting point is a result of W. Graustein \cite{Gr}.

\begin{theorem} \label{thm:E}
The average curvature of a plane oval is attained at least at four points.
\end{theorem}

 This implies the 4-vertex theorem: by Rolle's theorem, there is an extremum of the curvature between every two such points. 

Let $s$ be the arc length parameter on the positively oriented oval $\gamma$, $k(s)$ be the curvature, $R(s)=\frac{1}{k(s)}$ be the curvature radius, $L$ and $A$ be the length and the area bounded by the curve. Then the total curvature is $\int k(s)\, ds = 2\pi$, and the average curvature is $\bar k = \frac{2\pi}{L}$.

I shall present two proofs of Graustein's result in the Euclidean plane, starting with his own, and then discuss the extensions of this theorem to the spherical and hyperbolic geometries.
\medskip

{\bf Proof by Sturm-Hurwitz theorem} [Graustein]. Let $\alpha(s)$ be the direction of the tangent line of $\gamma$. Then $\gamma'(s)=(\cos\alpha(s),\sin\alpha(s))$ and $k(s)=\frac{d\alpha(s)}{ds}$ or, equivalently, $ds=R\,d\alpha$. One can use $\alpha$ to parameterize $\gamma$.

Since the curve is closed, 
$$
\int_0^L  (\cos\alpha(s),\sin\alpha(s))\ ds = 0 \ \ \ {\rm or}\ \ \ \int_0^{2\pi} R(\alpha) (\cos\alpha,\sin\alpha(s))\ d\alpha = 0.
$$
Hence the $2\pi$-periodic function $R(\alpha)$ is $L_2$-orthogonal to $\cos\alpha$ and $\sin\alpha$. 

Let $\bar R$ be the average value of $R(\alpha)$:
$$
\bar R = \frac{\int_0^{2\pi} R(\alpha)\ d\alpha}{2\pi}= \frac{\int_0^{L} dt}{2\pi} = \frac{L}{2\pi}.
$$
Then the function $R(\alpha) - \bar R$ is also $L_2$-orthogonal to constants. 

To conclude the proof, invoke the Sturm-Hurwitz theorem: {\it The number of sign changes of a $2\pi$-periodic function is not less than the number of sign changes of the first nontrivial harmonic in its Fourier expansion} (see, e.g., \cite{OT}, Appendix 8.1, for five different proofs). The Fourier expansion of the function $R(\alpha) - \bar R$ is free from the constant term and the first harmonics, hence this function has at least four zeros, the number of roots of the second harmonics. 

That is, $R$ assumes the value $\bar R=\frac{L}{2\pi}$ at least four times, and hence the curvature $k=\frac{1}{R}$ also assumes the value $\bar k = \frac{2\pi}{L}$ at least four times. 
\medskip

{\bf Proofs by wave propagation}. Consider the oval $\gamma$ as a source of light, and let $\gamma_t$ be the time-$t$ wavefront, 
the set of points reached by light at time $t$. We coorient $\g$ in the outward direction; then, for $t>0$, every point of $\g$ is moving with the unit speed along the outer normal ray to the curve, and for $t<0$, the points move along the normal lines inside the curve $\g$. The wavefronts $\g_t$ are equidistant from $\g=\g_0$ at distance $|t|$.

\begin{figure}[ht]
\centering
\includegraphics[width=.5\textwidth]{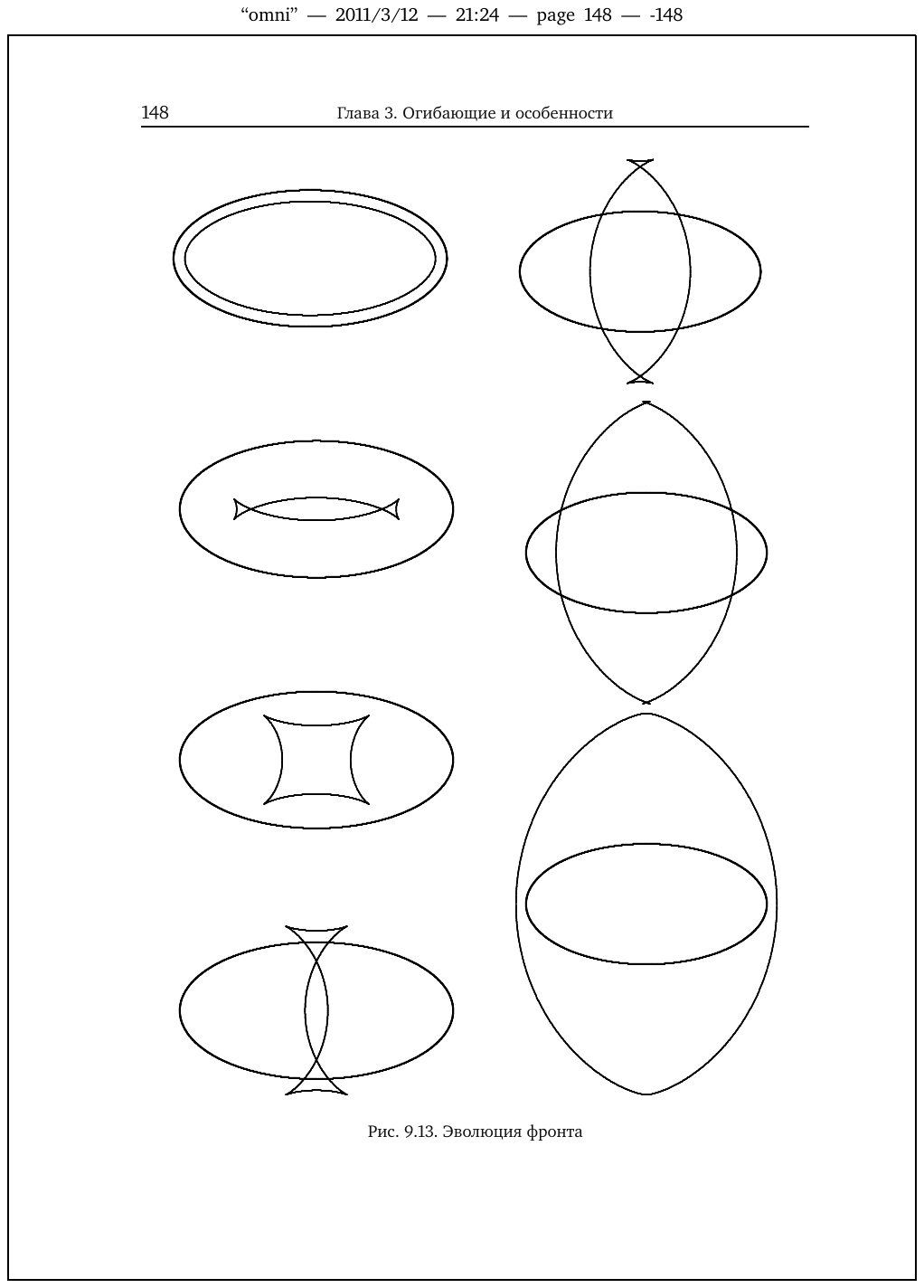}
\caption{Propagation of light inside an oval. Eventually, the curve inverts outside in.}
\label{invert}
\end{figure}

 For positive values of $t$, the curves $\g_t$ remain smooth and convex, but for negative $t$ they start to develop singularities, generically semicubical cusps, as can be seen in Figure \ref{invert} borrowed from \cite{FT}. 
 
We parameterize the equidistant curves $\g_t$ by the same parameter $s$:
$$
\g_t(s)=\g(s)+ t J\g'(s),
$$
where $J$ is the clockwise rotation by $90^\circ$ and hence $J\g'(s)$ is the unit outward normal.

\begin{figure}[ht] 
\centering
\includegraphics[width=.3\textwidth]{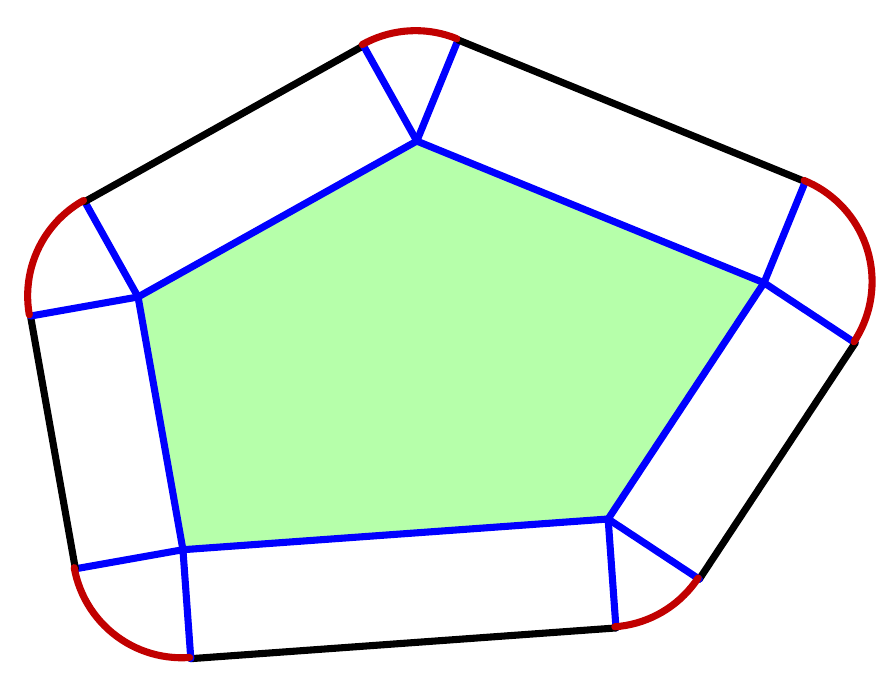}
\caption{Steiner's formulas for a polygon: the total length of the black segments is $L$, and the red arcs comprise a circle of radius $t$. }
\label{polyg}
\end{figure}

Let $k_t, R_t, L_t, A_t$ and $\bar k_t$ be the curvature, curvature radius, perimeter, area, and the average curvature of the curve $\g_t$. The dependences of these quantities on $t$ are as follows:
\begin{equation} \label{eq:ont}
R_t(s)=R(s)+t,\ \ k_t(s)=\frac{k(s)}{1+t k(s)},\ \ L_t=L+2\pi t,\ \ A_t=A+L t + \pi t^2.
\end{equation}
The first formula in (\ref{eq:ont}) follows from the Huygens principle, according to which the tangent elements of $\g$ move along the normals and remain tangent to the wavefronts $\g_t$. The second formula follows from the first one.
The  formulas for the length and area of $\g_t$ are classical and due to J. Steiner \cite{St}. For polygons, they are illustrated in Figure \ref{polyg}. A proof also follows from the obvious formulas
\begin{equation} \label{eq:St}
\frac{dL_t}{dt}= \int k_t(s)\ ds = 2\pi,\ \ \frac{dA_t}{dt} = L_t.
\end{equation}

For curves with cusps, the length is signed: the sign changes after traversing a cusps. Note also that the winding number of the wave fronts $\g_t$ (the number of turns of its tangent line) is 1 for all $t$.

\begin{remark} \label{rmk:iso1}
{\rm It is worth mentioning that the non-negative isoperimetric defect $L^2 - 4\pi A$ remains constant in the equidistant family:
$$
L_t^2-4\pi A_t = (L+2\pi t)^2 - 4\pi (A+L t + \pi t^2) = L^2 - 4\pi A.
$$
}
\end{remark}

The next lemma says that the points of an oval where the curvature equals the average curvature propagate in the equidistant family.

\begin{lemma} \label{lm:prop}
If $k(s)=\bar k$, then $k_t(s)=\bar k_t$ for all $t$.
\end{lemma}

\proof Note that 
$$
\bar k_t = \frac{2\pi}{L_t}=\frac{2\pi}{L+ 2\pi t}.
$$
Hence if
$k(s)=\frac{2\pi}{L}$, then
$$
k_t(s)=\frac{k(s)}{1+t k(s)} = \frac{2\pi}{L+ 2\pi t}=\bar k_t, 
$$
as needed.
\proofend

To prove Theorem \ref{thm:E}, consider the wavefront $\g_t$ with 
$$
t= - \frac{1}{\bar k} = - \frac{L}{2\pi}.
$$
The curvature at a point $\g(s)$ equals $\bar k$ if and only if the curvature of the respective point of $\g_t(s)$ is infinite, that is, this point is a cusp. Thus we need to show that $\g_t$ has at least four cusps. Two cusps are ``free": the average curvature must be attained at least twice. Working toward a contradiction, assume that there are only two cusps.

\begin{figure}[ht] 
\centering
\includegraphics[width=.3\textwidth]{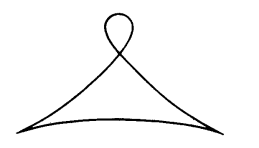}
\caption{This front cannot have zero length.}
\label{front}
\end{figure}

By formulas (\ref{eq:ont}), $L_t=0$. Note that the wave fronts are locally convex, that is, free from inflections: the curvature radius $R_t(s)=R(s)+t$ is never infinite. Thus the wave front $\g_t$ has the winding number 1 and does not have inflection points. Such a front is diffeomorphic to the one in Figure \ref{front} and, due to the triangle inequality, its length does not vanish. This is a contradiction.

Perhaps the claim that the curve in Figure \ref{front} cannot have zero signed length needs an explanation. It is consequence of the following well known result. 

\begin{lemma} \label{lm:conv}
Let $U$ be a domain with piecewise smooth boundary $\Gamma$, and $\gamma$ a closed convex curve strictly inside $U$.
Then the length of $\Gamma$ is greater than that of $\g$.
\end{lemma}

\begin{figure}[ht] 
\centering
\includegraphics[width=.4\textwidth]{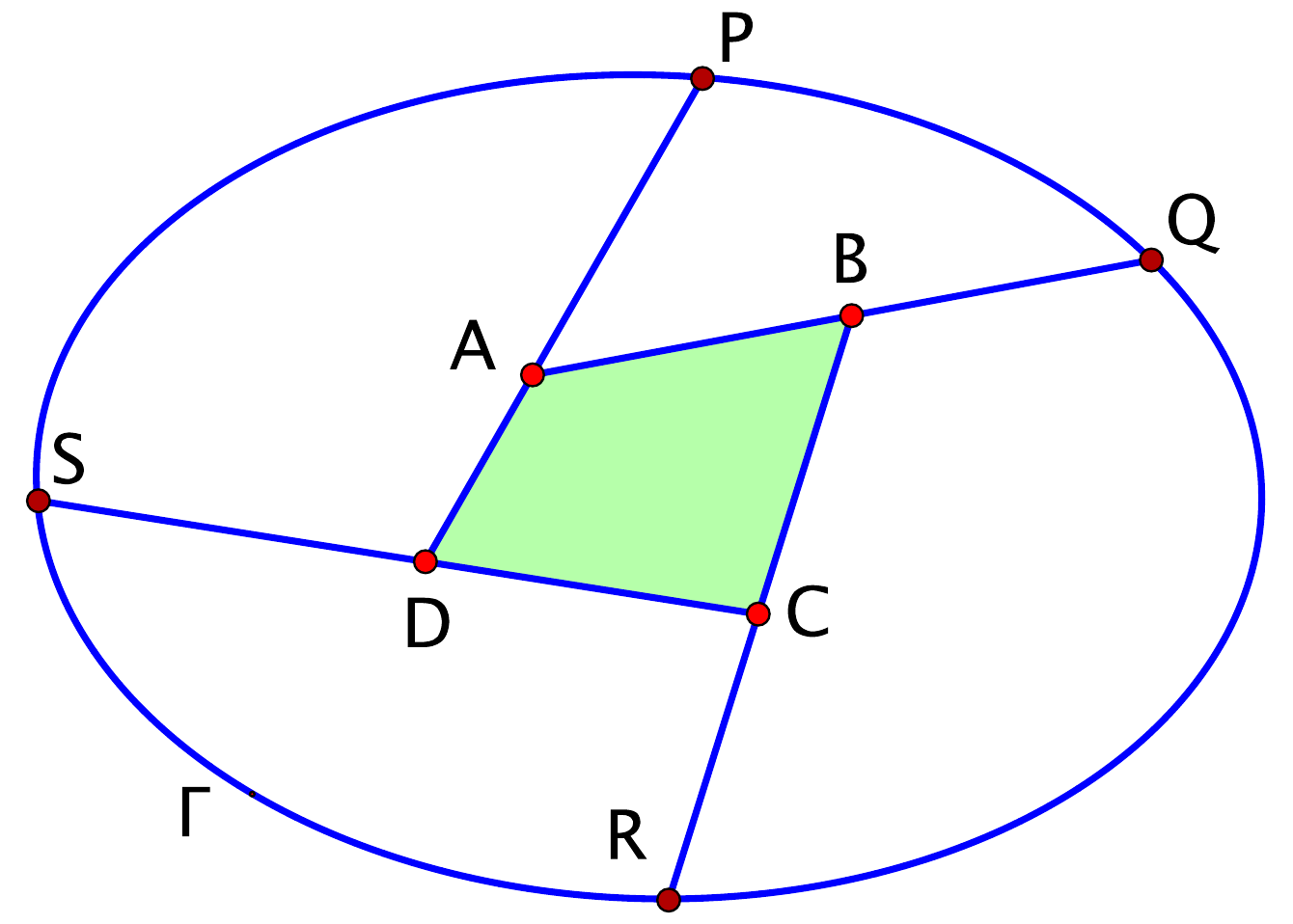}
\caption{The length of the curve $\Gamma$ is greater than the perimeter of the polygon $ABCD$.}
\label{tria}
\end{figure}

\proof One way of proving this claim is to approximate $\g$ by a convex polygon and use the triangle inequality, as illustrated in Figure \ref{tria}.
In this figure, one has
\begin{equation*}
\left\{\;
\begin{aligned}
|AP|+|\Gamma_{PQ}| > |AB|+|BQ|,\\
|BQ|+|\Gamma_{QR}| > |BC|+|CR|,\\
|QR|+|\Gamma_{RS}| > |CD|+|DS|,\\
|DS|+|\Gamma_{SP}| > |DA|+|AP|.
\end{aligned}\right.
\end{equation*}
Adding these inequalities and collecting terms yields $|\Gamma| > |\g|$.
\proofend

A more conceptual proof is by applying the Cauchy-Crofton formula; see, e.g.,  Lecture 19 of \cite{FT}. 

\begin{remark}
{\rm The curve shown in Figure \ref{front}  appeared in the literature on the 4-vertex theorem and related topics more than once: see \cite{Ar}, figure 30, and \cite{Ta1,Ta2,Va}.
}
\end{remark}

{\bf In the spherical geometry}. The wave propagation argument works in this case as well. We have a (geodesically) convex curve $\g(s)$ on the unit sphere, parameterized by the arc length and cooriented outwards. Such a curve lies in an open  hemisphere. Consider the family of equidistant curves $\g_t (s)$.

We need spherical versions of the equations (\ref{eq:ont}):
\begin{equation} \label{eq:sph}
\begin{aligned}
R_t(s)=R(s)+t&,\ \ k_t(s)=\cot R_t(s)=\cot (R(s)+t),\\ 
L_t=L\cos t + (2\pi-A)\sin t&,\ \ A_t=2\pi + L\sin t - (2\pi-A)\cos t.
\end{aligned}
\end{equation}
The second equation, relating the curvature to the curvature radius, is a standard formula of spherical geometry. The last two equations follow from the spherical analog of equations (\ref{eq:St}):
$$
\frac{dL_t}{dt}= \int k_t(s)\ ds = 2\pi -A,\ \ \frac{dA_t}{dt} = L_t,
$$
where the first equation is the Gauss-Bonnet formula. See \cite{An} for generalizations. 

Similarly to Remark \ref{rmk:iso1}, one has

\begin{remark} \label{rmk:iso2}
{\rm
The spherical isoperimetric inequality reads
$L^2 \geq A(4\pi-A).$
Like in the plane, the isoperimetric defect remains constant in the equidistant family:
$L_t^2-A_t(4\pi-A_t)$ is independent of $t$.
}
\end{remark}

According to the Gauss-Bonnet formula, $\bar k_t = \frac{2\pi-A_t}{L_t}$. Formulas (\ref{eq:sph}) imply an analog of Lemma \ref{lm:prop}: {\it If $k(s)=\bar k$, then $k_t(s)=\bar k_t$ for all $t$} (we omit this straightforward calculation).

Consider now the equdistant curve $\g_t$ with 
$$
t=\tan^{-1} \left(\frac{2\pi-A}{L}\right).
$$
Then $A_t=2\pi$, and $\bar k_t=0$. We need to show that this curve has at least four spherical inflection points.

We claim that the curve $\g_t$ is smooth and embedded. Indeed, if $\g_t$ has a cusp, then $k_t(s) = \infty$, and formulas (\ref{eq:sph}) imply that $\cot t = - k(s)$. However $\cot t = \frac{L}{2\pi-A} > 0$ because $\g$ is contained in a hemisphere. This is a contradiction.

\begin{figure}[ht] 
\centering
\includegraphics[width=.45\textwidth]{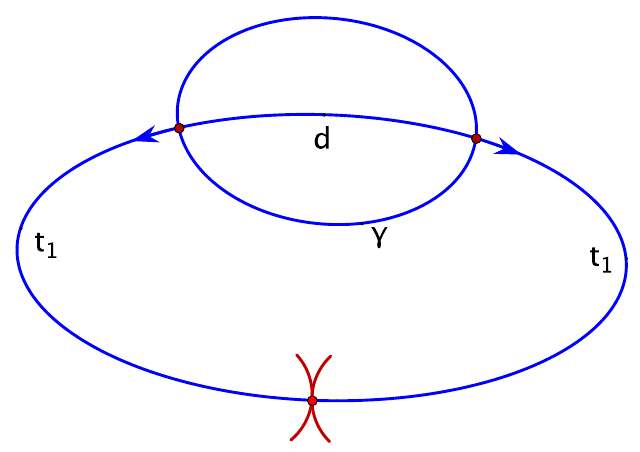}
\caption{One has $2t_1+d=2\pi$, hence $t_1=\pi-\frac{d}{2}> \frac{\pi}{2}$.}
\label{self}
\end{figure}

Likewise, if the curve $\g_t$ has self-intersections then, for some $t_1\le t$, the curve $\g_{t_1}$ developed a self-tangency. Note that $t_1\le t < \frac{\pi}{2}$. Consider Figure \ref{self}: $d$ is the spherical diameter of $\g$ (a geodesic chord perpendicular to the curve on both ends), and since $\g$ is contained in a hemisphere, $d<\pi$. This implies that $t_1>\frac{\pi}{2}$, a contradiction.

\begin{figure}[ht] 
\centering
\includegraphics[width=.2\textwidth]{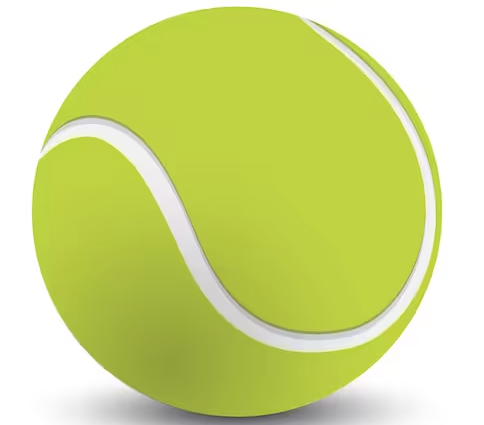}
\caption{Illustrating the ``tennis ball theorem".}
\label{tennis}
\end{figure}

Finally,  a coda: {\it a closed embedded spherical curve that bisects the area has at least four inflection points}. This is V. Arnold's ``tennis ball theorem", see \cite{Ar,Ar1}. The name is explained by Figure \ref{tennis}. We note, in passing, that this theorem follows from a stronger earlier result of B. Serge: {\it if a closed embedded spherical curve intersects every great circle, then it has at least four inflection points}, see \cite{Se}. 

\begin{remark} \label{rmk:torsion}
{\rm
It would be an omission not to mention the paper \cite{Ja} by S. Jackson devoted to a detailed study of the geometry of spherical curves, in particular, of their vertices. Let us mention but one result that has a similar flavor to out topic: 
the total torsion of a  closed spherical curve, considered as a space curve, vanishes, and if the curve is simple, then its torsion vanishes at  least four times (Section 6 of \cite{Ja}). The geodesic convexity of the curve is not needed.
}
\end{remark}

\bigskip
{\bf In the hyperbolic plane}. The formulas of spherical geometry have analogs in the hyperbolic plane: one replaces trigonometric functions by their hyperbolic counterparts and changes signs where needed. In particular, the hyperbolic analogs of formulas (\ref{eq:sph}) are as follows:
\begin{equation} \label{eq:hyp}
\begin{aligned}
R_t(s)=R(s)+t&,\ \ k_t(s)=\coth R_t(s)=\coth (R(s)+t),\\ 
L_t=L\cosh t + (2\pi+A)\sinh t&,\ \ A_t=-2\pi + L\sinh t + (2\pi+A)\cosh t.
\end{aligned}
\end{equation}
The last two formulas follow from the hyperbolic analog of equations (\ref{eq:St}):
$$
\frac{dL_t}{dt}= \int k_t(s)\ ds = 2\pi + A,\ \ \frac{dA_t}{dt} = L_t.
$$

We note that, as in Remarks \ref{rmk:iso1} and \ref{rmk:iso2}, the isoperimentric defect $L_t^2-A_t(4\pi+A_t)$ is independent of $t$. Note also that an analog of Lemma \ref{lm:prop} holds in the hyperbolic setting as well. 

This said,  formulas (\ref{eq:hyp}) work only for {\it horocyclically convex} curves, that is, the curves whose curvature is everywhere greater than 1. 

The reason is that, in the hyperbolic plane, a curve of a constant curvature $k$ is not necessarily a circle. Indeed, if $k>1$, it is a circle, if $k=1$, it is a horocycle, a circle of infinite radius tangent to the circle at infinity, and if $k<1$, this is an equdistant curve, the locus of points at a fixed distance from a geodesic. Accordingly, at every point a smooth curve  is approximated by a curve of constant curvature, but this is the osculating circle only if the curvature at this point is greater than 1. 

So assume that our curve $\g$ is horocyclically convex. Then the wave propagation proof works in the same way as in the Euclidean plane. Consider the projective (Beltrami-Cayley-Klein) model of $H^2$, represented by the unit disc, whose geodesics are its chords.

Consider the equidistant curve $\g_t$ with 
$$
t=- \coth^{-1} \left(\frac{2\pi+A}{L}\right),
$$
implying that $L_t=0$. 

The curve $\g_t$ is regularly homotopic to $\g$, so its winding number is 1. It is also free from inflections -- this is where we use the horocyclic convexity of $\g$. Therefore it is diffeomorphic to the wavefront in Figure \ref{front} and cannot have zero signed length.
\medskip

{\bf And a counter-example}. The next example, suggested by A. Akopyan, shows that one cannot relax the assumption of horocyclic convexity. 

\begin{figure}[ht] 
\centering
\includegraphics[width=.4\textwidth]{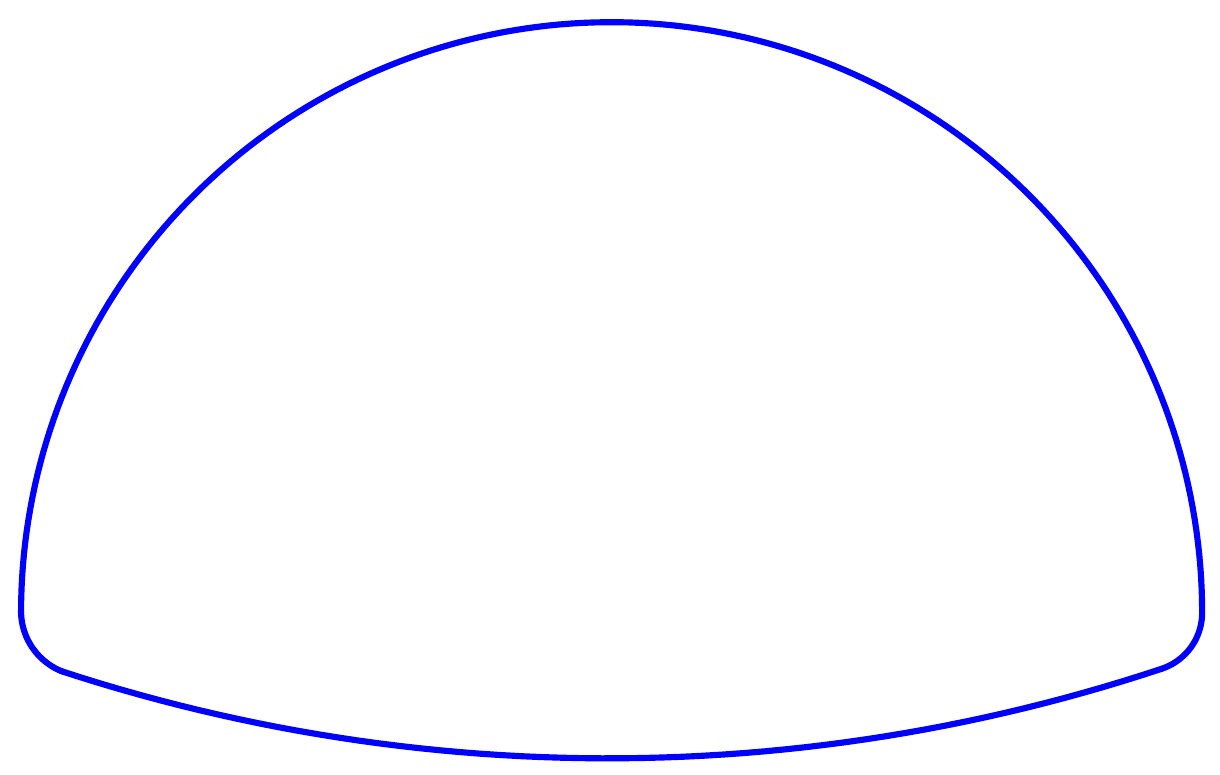}
\caption{A rounded semicircle.}
\label{example}
\end{figure}

Consider a semicircle of radius $r$ in the hyperbolic plane and smooth it by replacing the corners by curves of a very large curvature and the diameter by a curve of a very small curvature, see Figure \ref{example}.

 The average curvature of the semicircle is calculated using formulas (\ref{eq:hyp}):
$$
\bar k = \frac{2\pi+A}{L} = \frac{\pi(1+\cosh r)}{2r+\pi \sinh r},
$$
where $L$ and $A$  are the length of the semicircle and the area bounded by it.
At the same time, the curvature of the circle is $\coth r$.

If $r$ is sufficiently large, then $\bar k < \coth r$ (it is enough to have $r>1.386$). If the approximation of the semicircle is fine enough,
then the average curvature $\bar k$ is attained on the nearly flat arc replacing the diameter, and  one can arrange it to happen only twice. 
\medskip

Let me conclude with a question: {\it do the above 4-point results hold in the Euclidean and spherical geometries for simple closed, but not necessarily convex, curves?}

\bigskip
{\bf Acknowledgements}. I am grateful to A. Akopyan, P. Albers, V. Ovsienko, and R. Schwartz  for stimulating discussions.
I was supported by NSF grant DMS-2404535, and by a Mercator fellowship within the Transregional Colloborative Research Center CRC/TRR 191 (281071066). It is a pleasure to acknowledge the hospitality of the Heidelberg University.

\end{document}